\documentclass[12pt]{amsart}
\usepackage{amsmath,amsfonts,amssymb}

\def\ssum{\mathop{\sum\!\sum}}

\def\ssssum{\mathop{\ssum\!\ssum}}

\def\le{\leqslant}

\def\ge{\geqslant}

\def\ab{a^6+b^2}

\def\vp{\varphi}
\def\frak{\mathfrak}

\renewcommand{\mod}{\mathop{\rm{mod}}}
\def\ab{\mathbb{A}}
\def\bb{\mathbb{B}}

\def\cA{\mathcal{A}}

\def\cC{\mathcal{C}}
\def\cD{\mathcal{D}}
\def\cP{\mathcal{P}}
\def\cW{\mathcal{W}}

\numberwithin{equation}{section}

\theoremstyle{definition}

\begin{document}
\title{\bf Sifting for small primes \\ from an arithmetic progression}
%\title{\bf Sifting the least prime \\ from an arithmetic progression}
%\title{\bf Sifting the least prime}

\author[Friedlander]{J.B. Friedlander$^*$}
\thanks{$^*$\ Supported in part by NSERC grant A5123}
\author[Iwaniec]{H. Iwaniec}

\maketitle

\dedicatory \quad\quad\quad\quad{\sl On the fiftieth anniversary of Chen's Goldbach theorem}

\medskip

{\bf Abstract:} In this work and its sister paper [FI3] we give a new proof of the famous Linnik theorem bounding the least prime in an arithmetic progression. Using sieve machinery in both papers, we are able to dispense with the log-free zero density bounds and the repulsion property of exceptional zeros, two deep innovations begun by Linnik and relied on in earlier proofs.
\footnote{key words: sieve, small primes, Linnik theorem} 

\section{\bf Introduction: A bit of history}
The theorem of Yu.V. Linnik [L] asserts that, for integers $q\ge 2$, $(a,q)=1$, the least prime $p \equiv a(\mod q)$ satisfies 
\begin{equation}\label{eq:1.1}
p_{\rm min}(q,a) \ll q^L ,
\end{equation}
where $L$ and the implied constant are absolute. 

The Riemann Hypothesis for Dirichlet $L$-functions $L(s,\chi)$, $\chi (\mod q)$ yields 
\begin{equation}\label{eq:1.2}
p_{\rm min}(q,a) \ll (q\log q)^2 .
\end{equation} 
Linnik's unconditional arguments make use of the Deuring-Heilbronn phenomenon which describes a repulsion effect of the exceptional zero on all of the other zeros. Even if the exceptional zero does not exist, the standard arguments, then and now, do not produce ~\eqref{eq:1.1} with an absolute constant $L$. To this end, Linnik established the so-called log-free density bound 
\begin{equation}\label{eq:1.3}
\sum_{\chi (\mod q)}\#\{\rho =\beta + i\gamma; L(\rho,\chi)=0,\, \alpha \le \beta <1, |\gamma|\le q\}\ll q^{c(1-\alpha)} ,
\end{equation}
where $c$ is an absolute constant. 
Moreover, in the presence of the exceptional zero, say $\beta_1$, with $\eta_1= (1-\beta_1)\log q$ being sufficiently small, Linnik was able to expand the classical zero-free region, with the result that all the other zeros $\rho=\beta + i\gamma$ with $|\gamma| \le q$ satisfy 
\begin{equation}\label{eq:1.4}
(1-\beta)\log q \gg \log 1/\eta_1 
\end{equation}
with an absolute implied constant. This is a quantitative form of the Deuring-Heilbronn phenomenon. Later, Bombieri [B] combined these two crucial estimations of Linnik into a single statement, namely that the number of all zeros $\rho = \beta + i\gamma$ with $\alpha \le \beta <1$, $|\gamma| \le q$ of all $L(s,\chi)$ other than $\beta_1$ satisfies 
\begin{equation}\label{eq:1.5}
\sum_{\chi (\mod q)}N^*(\alpha, \chi, q) \le c_1(1-\beta_1)(\log q) q^{c(1-\alpha)} ,
\end{equation}
where the constants $c_1$, $c$ are absolute. See Th\'eor\`eme 14 of [B] for somewhat stronger statements. Bombieri employs the Turan power sum method, essentially in place of the sieve. 
 
Sieve ideas in conjunction with the mollification of $L$-functions constitute the driving force of Linnik's work. His arguments have been refined substantially over time with new innovative elements (particularly those of Heath-Brown), giving concrete values of $L$. Here are some remarkable records:

\[
\begin{array}{cccccc}
L & = & 10^5  & ({\rm Pan}, 1957)\\
L & = & 80    & ({\rm Jutila}, 1977)\\
L & = & 20    & ({\rm Graham}, 1981)\\
L & = & 13.5  & ({\rm Chen\, and \, Liu}, 1989)\\
L & = &  5.5  & ({\rm Heath-Brown}, 1992)\\
L & = &  5.18 & ({\rm Xylouris}, 2009) 
\end{array} 
\]

It is interesting that, if permitted the assistance of the exceptional character $\chi_1 (\mod q)$, one can accomplish a lot more for this problem than with the Riemann Hypothesis. For example, we obtained~\eqref{eq:1.1} with $L= 2 -1/55$ under suitable conditions, see [FI1]. If one tries to count primes $p \equiv a (\mod q)$ directly using the Riemann Hypothesis, by means of the explicit formula 
\begin{equation*}
\begin{aligned}
\psi_f(x;q,a) & = \sum_{n\equiv a (\mod q)}f(n/x)\Lambda(n) \\
& = {\hat f}(1)\frac{x}{\varphi (q)} -
\frac{1}{\varphi (q)}\sum_{\chi \neq \chi_0}\bar \chi(a)\sum_{\rho}{\hat f}(\rho)x^{\rho} \\
& = {\hat f}(1)\frac{x}{\varphi (q)} +O(\sqrt x \log x) 
\end{aligned}
\end{equation*}
with $f$ smooth and compactly supported on $(0, \infty)$, 
then one may get~\eqref{eq:1.1} with $L< 2$ provided that there is a considerable cancellation in the terms $x^{\rho}$ as $\rho$ runs over low-lying zeros $\rho = \beta + i\gamma$.
% In 1976, H.L. Montgomery conjectured that 
%\begin{equation*}
%\psi(x;q,a)  = \sum_{\substack{n\le x \\ n\equiv a (\mod q)}}\Lambda(n) = 
%\frac{x}{\varphi (q)} +O(q^{-\frac 12}x^{\frac 12 +\varepsilon})
%\end{equation*}
% which would yield
One actually expects that~\eqref{eq:1.1} holds with any $L>1$. The distribution of primes $p \equiv a (\mod q)$ in the vicinity of $q$ is tricky; cf. [FG]. 

In these notes we do not attempt to obtain a reasonably small Linnik constant $L$. In fact, at the end of the day we arrive at~\eqref{eq:1.1} with the huge value $L=75,744,000$. So, what motivates us to produce this long proof of Linnik's bound? Our first reason is to demonstrate how the sieve machinery can succeed in this problem when fueled by neither the log-free density estimation nor the repulsion property of the exceptional zero. However, we do appeal to the classical zero-free region, specifically Lemma A1 with $\eta = 1/657$.
 
Moreover, these notes provide more details to our arguments in Chapter 24 of [FI2]. To diversify a bit we treat the case of the exceptional zero by Selberg's lower-bound sieve rather than by a combinatorial sieve.  Of course, the results are essentially the same, but we get nicer expressions. We have put the treatment of that case into a separate note [FI3]. 

Our second incentive to produce these notes is the desire to fix a defect in the technical arrangements of the partial sums $\psi_{\frak{X}}(x;q,a)$ in [FI2], a correction of which we have no other opportunity to publish due to an AMS copyright issue. The required correction is embarrassing but not difficult. The subdivision there of $n=pp_1p_2p_3p_4 \equiv a (\mod q)$ should have been made over the segments

\begin{equation}\label{eq:1.6}
D^{\frac16} < p_1,p_2 < D^{\frac15}, \quad x_3<p_3\le 2x_3, \quad x_4<p_4\le 2x_4
\end{equation}
rather than over $x_j<p_j \le 2x_j$  for all $j=1, 2, 3, 4$. That is, we should not have put $p_1, p_2$ into dyadic segments as we did for $p_3, p_4$. However, adjusting the existing arguments to the larger segments~\eqref{eq:1.6} is simple with obvious modifications. The key difference appears in the factor $h(p/x)$ of (24.47) in [FI2] which needs to be changed to $h(\log p/\log x)$. In other words, the first primes $p_1, p_2$ should have been logarithmically scaled. See how we proceed in Subsection 2.1 of these notes. 

\section{\bf The linear beta-sieve}

Let $\cA = (a_n)$ be a finite sequence of real numbers $a_n\ge 0$, supported on $x<n\le 2x$, $x$ large. Our goal is to estimate 

\begin{equation}\label{eq:2.1} 
S(\cA ) =\sum_p a_p .
\end{equation}
It is well understood that, due to the parity barrier, one cannot obtain a positive lower bound for $S(\cA )$ within the basic sieve methods. In this section we reduce the problem for primes to that for products of five primes, that is to the estimation of  
\begin{equation}\label{eq:2.2} 
S_5(\cA ) =\sum_p\sum_{p_1}\sum_{p_2}\sum_{p_3}\sum_{p_4} a_{pp_1p_2p_3p_4} ,
\end{equation}
where the $p_j$ run over specific segments, short in the logarithmic scale. We do not lose the correct order of magnitude. Note that $n=pp_1p_2p_3p_4$, which we call a prime quintet, consists of integers with an odd number of prime factors. 

Let $P(z)$ denote the product of primes $p<z$ in some set $\cP$. Sieve methods are capable of producing upper and lower bounds for the sifting function
\begin{equation}\label{eq:2.3} 
S(\cA,z ) =\sum_{(n,P(z))=1} a_n ,
\end{equation}
but only in a limited range.  The basic assumption needed for sifting is that the congruence sums
\begin{equation}\label{eq:2.4} 
A_d =\sum_{n\equiv 0(\mod d)} a_n 
\end{equation}
can, on average over $d$, be well-approximated by a simple model $g(d)X$, where $g(d)$ is multiplicative with $0\le g(p)<1$ and $X$ is fixed. This requirement means that the error terms in 
\begin{equation}\label{eq:2.5} 
A_d =g(d)X +r_d 
\end{equation}
are relatively small. In particular, we must choose $X\sim A_1$. The probability of seeing $a_n$ with $(n,P(z))=1$ is given by the product 
\begin{equation}\label{eq:2.6} 
V(z) =\prod_{p|P(z)} \bigl(1-g(p)\bigr)
\end{equation}
and sieve methods produce bounds of the true order of magnitude 
\begin{equation}\label{eq:2.7} 
S(\cA,z) \asymp XV(z) 
\end{equation}
provided that we can control the error terms $r_d$. In practice, we can do so only for moduli up to some level $y\le x$. We need, say, 
\begin{equation}\label{eq:2.8} 
R(y)= \sum_{\substack{d<y\\d|P(z)}}|r_d| \le XV(z)(\log y)^{-1} .
\end{equation}

Let $\cA_d$ be the subsequence of $\cA$ of elements $a_n$ with $n\equiv 0 (\mod d)$. We start from the combinatorial identity
\begin{equation}\label{eq:2.9} 
S(\cA,z)= S^-(\cA,z) + S_2(\cA,z) + S_4(\cA,z) +\ldots  
\end{equation}
which is derived by a process of inclusion-exclusion; see (12.10) and (12.11) for the linear beta-sieve in [FI2]. This holds for $z=\sqrt y$ with
\begin{equation}\label{eq:2.10} 
S^-(\cA,z) =\sum_{d|P(z)}\lambda^-_dA_d 
\end{equation}
where $\lambda^-_d$ are the sieve weights. We have
\begin{equation}\label{eq:2.11} 
S_4(\cA,z) = \ssssum_{p_4<p_3<p_2<p_1}S(\cA_{p_1p_2p_3p_4},p_4) 
\end{equation}
where the prime variables of summation are restricted by two inequalities
\begin{equation}\label{eq:2.12} 
p_1p_2^3 <y, \quad p_1p_2p_3p_4^3\ge y .
\end{equation}
The other $S_n(\cA,z)$ are defined similarly but we are going to drop them due to positivity. Actually, we also discard some parts of~\eqref{eq:2.12}. Specifically we take
\begin{equation}\label{eq:2.13} 
x^{\frac16} <p_j < x^{\frac15}, \quad j=1,2,3,4
\end{equation}
and $a_n$ with $n=pp_1p_2p_3p_4$, $p$ prime only. We assume $x^{\frac45}<y<x$, so that the restrictions~\eqref{eq:2.12} hold automatically because $x\le n\le 2x$.
Observe that $p$ is the largest prime in the prime quintet $n=pp_1p_2p_3p_4$; specifically we have 
\begin{equation}\label{eq:2.14} 
x^{\frac15} <p < 2x^{\frac13} .
\end{equation}

Let $Q(\cA)$ be the following specific choice for~\eqref{eq:2.2}
\begin{equation}\label{eq:2.15} 
Q(\cA) =\sum_p\sum_{p_1}\sum_{p_2}\sum_{p_3}\sum_{p_4} a_{pp_1p_2p_3p_4} 
\end{equation}
where all prime variables are distinct and $p_j$ run freely over the segment~\eqref{eq:2.13} without ordering. The variable $p$ is restricted only by the support of $\cA=(a_n)$ because~\eqref{eq:2.14} is redundant. We obtain the following combinatorial inequality: 
\begin{equation}\label{eq:2.16} 
S(\cA,z) \ge S^-(\cA,z) +\frac{1}{24}Q(\cA) , 
\quad {\rm if}\,\, y^{1/4}\le z\le y^{1/2} .
\end{equation}

The linear beta-sieve in Chapter 12 of [FI2] requires the density function $g(d)$ to satisfy the upper bound 
\begin{equation}\label{eq:2.17} 
\prod_{w\le p <z}\bigl(1-g(p)\bigr)^{-1} \le\frac{\log z}{\log w}\bigl(1+\frac{\ell}{\log w}\bigr)
\end{equation}
for every $2\le w <z$ with some constant $\ell$. In our case, we shall verify~\eqref{eq:2.17} with an absolute constant $\ell$. Then we have 
 \begin{equation}\label{eq:2.18} 
S^-(\cA,z)\ge \bigl(f_1(s)+O\bigl((\log y)^{-1/6}\bigr)\bigr)XV(z) - R(y)
\end{equation}
where $s=(\log y)/\log z$ and $f_1(s) =2e^{\gamma}s^{-1}\log(s-1)$, if $2\le s \le 4$. We take $z={\sqrt y}$, so $s=2$, $f_1(2)=0$ and~\eqref{eq:2.18} gives nothing positive. Therefore, we are left with the following lower bound
\begin{equation}\label{eq:2.19} 
S(\cA,\sqrt y)\ge \frac{1}{24}Q(\cA) + O\bigl(XV(\sqrt y)(\log y)^{-1/6}\bigr)\bigr).
\end{equation}

On the other hand, $S(\cA)=S(\cA,\sqrt {2x})$ is close to $S(\cA,{\sqrt y})$ if the admissible level of distribution $y$ is close to $x$ in the logarithmic scale. Put
\begin{equation}\label{eq:2.20} 
y=x^{\theta}, \quad \frac45 < \theta <1 .
\end{equation}
The difference
\begin{equation}\label{eq:2.21} 
S(\cA,\sqrt y) -S(\cA,{\sqrt 2x})=\sum_{\sqrt y \le p < \sqrt{2x}}S(\cA_p,p) ,
\end{equation}
counting $a_n$ over prime duos $pp'$ can be estimated successfully by any upper-bound sieve method. Our sieve conditions~\eqref{eq:2.8} and~\eqref{eq:2.17} allow us to apply (12.12) of [FI2] for every subsequence $\cA_p$. We replace $X$ in~\eqref{eq:2.5} by $g(p)X$, $r_d$ by $r_{pd}$ and choose the level of distribution for $\cA_p$ to be $y/p$. We obtain $S(\cA_p,p) \le S(\cA_p,y/p)$ and 
\begin{equation*}
S(\cA_p,y/p)\le g(p)V(y/p)X\bigl\{F_1(1) + O\bigl((\log y)^{-1/6}\bigr)\bigr\} +\sum|r_{pd}| ,
\end{equation*}
where $F_1(1)= 2 e^{\gamma}$ and the sum of the error terms $r_{pd}$ runs over $d|P(p)$, $d\le y/p$. By~\eqref{eq:2.17} we have
\begin{equation*}
V(y/p)\le V(y)\frac{\log y}{\log(y/p)}\bigl(1+ O\bigl(\frac{1}{\log y}\bigr)\bigl) . 
\end{equation*}
Next, we apply (5.47) of [FI2] as follows: 
\begin{equation*}
\begin{aligned}
\sum_{\sqrt y < p < \sqrt {2x}} g(p)(\log y/p)^{-1} & < \int_{\sqrt y}^{\sqrt {2x}}(\log y/w)^{-1}\frac{d\log w}{\log w} +O\bigl( (\log y)^{-2}\bigr) \\
& = \bigl( \log \frac{1}{2\theta -1} + O(1/\log y)\bigr)(\log y)^{-1} .
\end{aligned}
\end{equation*}
Moreover, we write by Mertens' formula, 
\begin{equation*}
V(y)\log y =e^{-\gamma}H\bigl(1+O(1/\log y)\bigr)
\end{equation*}
where 
\begin{equation}\label{eq:2.22} 
H= H(y) = \prod_{p <y}\bigl(1- g(p)\bigr)\bigl( 1-\frac{1}{p}\bigr)^{-1}. 
\end{equation}
Inputting the above estimates into~\eqref{eq:2.21} and~\eqref{eq:2.19}, we obtain the following result.

\smallskip 
%\begin{theorem}\label{thm:2.1} 
{\bf Theorem 2.1:}
{\sl Suppose~\eqref{eq:2.17} holds for every $2\le w < z \le\sqrt{2x}$ with an absolute constant $\ell$ and the remainder $R(y)$ of level~\eqref{eq:2.20} satisfies~\eqref{eq:2.8}. Then we have the following inequality.
\begin{equation}\label{eq:2.23} 
S(\cA)\ge \frac{1}{24}Q(\cA) -\frac{2HX}{\log y} \log \frac{1}{2\theta -1}
+ O\bigl(X(\log y)^{-7/6}\bigr), 
\end{equation}
where $Q(\cA)$ is given by~\eqref{eq:2.15} and the implied constant is absolute.} 
%\end{theorem} 
\smallskip

{\bf Remarks:} We shall succeed to show that $Q(\cA)\ge cHX(\log y)^{-1}$ with an absolute constant $c>0$, but it will be a small constant because we took account of only a small portion of the combinatorial decomposition~\eqref{eq:2.9}. Nevertheless, we still obtain a positive lower bound for $S(\cA)$ because $\theta$ will be very close to $1$. 

\smallskip

From now on we are interested in an arithmetic progression. Let $q$ be large, $(a,q)=1$ and $x\ge q^6$. We take $\cA = (a_n)$ with 
\begin{equation}\label{eq:2.24} 
a_n = f(n/x) , \quad {\rm if} \,\, n\equiv a (\mod q)
\end{equation}
and $a_n = 0$ otherwise. 
Here, $f(u) \ge 0$ is a smooth function supported on $1\le u\le 2$ with
\begin{equation}\label{eq:2.25} 
{\hat f}(0)=\int f(u)du >0 .
\end{equation}
Our sequence $\cA=(a_n)$ satisfies the sieve conditions with 
\begin{equation}\label{eq:2.26} 
X={\hat f}(0)x/q 
\end{equation}
and density function 
\begin{equation}\label{eq:2.27} 
g(d)=d^{-1} \,\, {\rm if} \,\, (d,q)=1 , \quad g(d)=0 \,\, {\rm if}\,\, (d,q)\neq 1 . 
\end{equation}
Hence, the product~\eqref{eq:2.22} is constant, namely $H= q/\varphi (q)$ because $y\ge q$. Now, $S(\cA)$ becomes 
\begin{equation}\label{eq:2.28} 
\pi_f(x; q,a) =\sum_{p\equiv a (\mod q)} f(p/x) 
\end{equation}
and~\eqref{eq:2.23} yields
\begin{equation}\label{eq:2.29} 
\pi_f(x; q,a) \ge \frac{1}{24}Q(\cA) - {\hat f} (0)\bigl(\log\frac{1}{2\theta -1}\bigr)\frac{2x}{\varphi(q)\log y}
+ O\bigl(\frac{x}{q}(\log x)^{-7/6}\bigr) .
\end{equation}

The error terms are absolutely bounded and $R(y)\ll y$ so we choose 
\begin{equation}\label{eq:2.30} 
y=x/q(\log x)^3 = x^{\theta} 
\end{equation}
and check that the remainder satisfies~\eqref{eq:2.8}. 

Now $Q(\cA)$ in~\eqref{eq:2.23} for our sequence~\eqref{eq:2.24} reads as 
\begin{equation}\label{eq:2.31} 
Q(\cA) = \sum_p\ssssum_{\substack{pp_1p_2p_3p_4\equiv a (\mod q)\\x^{1/6}< p_j <x^{1/5},\,\, j=1,2,3,4}
}f(pp_1p_2p_3p_4/x) .
\end{equation}

At this point we can predict the limit of our output. Definitely, we cannot do better than the asymptotic formula for $\vp(q)Q(\cA)$:
\begin{equation*}
\begin{aligned}
%\vp(q)Q(\cA) & \sim 
& \ssssum\int f(vp_1p_2p_3p_4/x)(\log v)^{-1}dv \\
& \sim  {\hat f}(0)x\ssssum (p_1p_2p_3p_4)^{-1}\bigl(\log (x/p_1p_2p_3p_4)\bigr)^{-1} \\
& \sim {\hat f}(0)\frac{x}{\log x}{\iiiint} 
(1-u_1-u_2-u_3-u_4)^{-1}\frac{du_1}{u_1}\frac{du_2}{u_2}\frac{du_3}{u_3}\frac{du_4}{u_4} \\
& <{\hat f}(0)\frac{5x}{\log x}\bigl(\log \frac65\bigr)^4 , 
\end{aligned}
\end{equation*}
where ${\frac{\theta}{6} <u_j < \frac{\theta}{5}}$ ,
Hence, if the lower bound in~\eqref{eq:2.29} is positive, then 
\begin{equation*}
\frac{5}{24}\bigl(\log\frac{6}{5}\bigr)^4 > 2\log \frac{1}{2\theta -1}, \quad  
{\rm so}\,\, 1-\theta < 6\cdot 10^{-5} .
\end{equation*}
This shows that we cannot produce $p\equiv a (\mod q)$ smaller than $q^{10^5/6}$. 

\medskip 

%\center {2.1 \, {\bf Technical preparations:}}
2.1 \, {\bf Technical preparations:}

\smallskip
   
Before detecting the congruence $n=pp_1p_2p_3p_4\equiv a (\mod q)$ by characters, we take advantage of positivity to separate the variables $p_3$, $p_4$ from the crop function $f(n/x)$. To this end we subdivide $p_3$, $p_4$ into short segments 
\begin{equation}\label{eq:2.32} 
C<p_3\le \lambda C, \quad D<p_4 \le \lambda D
\end{equation}
with 
\begin{equation}\label{eq:2.33} 
P\le C, D \le\lambda^{-1} P^{6/5}, \quad P=x^{1/6}. 
\end{equation}
The constant $\lambda$ is just slightly larger than $1$. The complete interval~\eqref{eq:2.13} is  covered by the union of the segments~\eqref{eq:2.32} with
\begin{equation*}
C=\lambda^mP, \, D=\lambda^nP, \,\, 0\le m, n \le (\log P)/5\log \lambda . 
%C=\lambda^mP$, $D=\lambda^nP$, $ 0\le m, n \le (\log P)/5\log \lambda . 
\end{equation*}
 
%\smallskip

We want to replace $f(pp_1p_2p_3p_4/x)$ by $f(pp_1p_2CD/x)$. To this end we introduce the function $g(u)=uf'(u)$ and apply the following: 
\begin{equation}\label{eq:2.34} 
\begin{aligned}
f(pp_1p_2p_3p_4/x) & = f(pp_1p_2CD/x) +\int_1^{p_3p_4/CD}g(pp_1p_2CDt/x)t^{-1}dt \\
& \ge  f(pp_1p_2CD/x) -\int_1^{\lambda^2}|g(pp_1p_2CDt/x)|t^{-1}dt \\
& = f(pp_1p_2CD/x) - h(pp_1p_2CD/x) ,
\end{aligned}
\end{equation}
where 
\begin{equation}\label{eq:2.35} 
h(u) = \int_1^{\lambda^2} |g(ut)|t^{-1}dt =u\int_1^{\lambda^2} |f'(ut)|dt\ll \log \lambda . 
\end{equation}
Note that $h(u)$ is supported on $\lambda^{-2}\le u\le 2$  and is small if $\lambda$ is very close to $1$. 

We put 
\begin{equation}\label{eq:2.36} 
Q_f(C,D) = \sum_p\ssssum_{\substack{pp_1p_2p_3p_4\equiv a (\mod q) \\ P < p_1, p_2 <P^{6/5} \\ C<p_3\le \lambda C,\, D<p_4\le \lambda D}}
f(pp_1p_2CD/x) .
\end{equation}
Similarly, we define $Q_h(C,D)$ as in~\eqref{eq:2.36} with $f(u)$ replaced by $h(u)$. Summing over $C=\lambda^mP$, $D=\lambda^nP$, in the segment~\eqref{eq:2.33}, we get
\begin{equation}\label{eq:2.37} 
Q(\cA)\ge \sum_C\sum_D\bigl(Q_f(c,D) - Q_h(C,D)\bigr) .
\end{equation}

The case of $Q_h(C,D)$ can be treated as $Q_f(C,D)$, but we can also estimate $Q_h(C,D)$ directly by applying the Brun-Titchmarsh theorem in the variable $p$. Let $\tau$ be a quantity $O(1/\log P)$, not always the same one. We find

\begin{equation*}
\sum_{p\equiv\alpha (\mod q)}h(pp_1p_2CD/x)< {\hat h}(0)\frac{(2+\tau) x}{\vp (q)}(\log x/qp_1p_2CD)^{-1}
\end{equation*}
where $qpp_1p_2CD<qP^{24/5}=qx^{4/5}\le x^{29/30}$. Hence
\begin{equation}\label{eq:2.38} 
Q_h(C,D)<{\hat h}(0)\frac{(2+\tau) }{\vp (q)}(\log \tfrac65)^2\bigl(\frac{\log\lambda}{\log P}\bigr)^2\frac{30x}{\log x}
\end{equation} 
and the contribution of $Q_h(C,D)$ to~\eqref{eq:2.37} is $\ll {\hat h}(0)x/\vp(q)\log x$, which is insignificant because ${\hat h}(0)=\int h(u)du= (1-\lambda^{-2})\int u|f'(u)|du$ is small for $\lambda$ close to $1$. Actually, this contribution will be discarded entirely after rounding up (or down) the constants in the other contributions. 

\section{\bf Dual sums: primes versus zeros}

Let $\rho = \beta + i\gamma$ run over complex numbers with $0\le \beta \le 1$ and $|\gamma| \le T$. For $T\ge 1$ and $P\ge 3$ we put 
\begin{equation}\label{eq:3.1} 
V=\max_{|\gamma|\le T}\sum_{\rho}\bigl(1+ (1-\beta)\log P\bigr)^{-1}\bigl(1+(\gamma -t)^2(\log P)^2\bigr)^{-1} .
\end{equation} 

\smallskip

{\bf Lemma 3.1:} {\sl For arbitrary complex numbers $a_p$ supported on primes in the segment $P\le p\le P^{6/5}$ we have 
\begin{equation}\label{eq:3.2} 
\sum_{\rho}P^{{\frac 52}(\beta -1)}\big|\sum_pa_p p^{-\rho}\big|^2 \le 
\bigl(1387V+O(NT(\log P)^{-4})\bigr)\sum_p|a_p|^2p^{-1}
\end{equation} 
where $N=N(T)$ denotes the number of $\rho$'s counted with multiplicity and the implied constant is absolute. }

\smallskip 

{\bf Proof:} The assertion~\eqref{eq:3.2} is, by duality, equivalent to the statement that the inequality 
\begin{equation}\label{eq:3.3} 
\sum_{P\le p\le P^{6/5}}\big|\sum_{\rho}z_{\rho} p^{{\frac12}-\rho}P^{{\frac54}(\beta -1)}\bigr|^2 \le 
\bigl(1387V+O(NT(\log P)^{-4})\bigr)\sum_{\rho}|z_{\rho}|^2
\end{equation} 
holds true for all complex numbers $z_{\rho}$. 

For the proof of the dual inequality~\eqref{eq:3.3} we majorize the summation over $p$   by attaching the weights $k(\log n/\log P)\Lambda (n)/\log P$ where $k(u)\ge 0$ is a twice differentiable function supported on $\frac 14 \le u \le \frac 54$ with $k(u) \ge 1$, if $1\le u\le \frac 65$. Then, the left side of~\eqref{eq:3.3} is bounded by
\begin{equation}\label{eq:3.4} 
\sum_{\rho}\sum_{\rho'}|z_{\rho}z_{\rho'}||K(s)|P^{{\frac 54}(\sigma -2)}
\end{equation} 
where $s=\bar\rho +\rho' = \sigma + it$ with $\sigma =\beta +\beta' \le 2$, $t=\gamma' -\gamma$, $|t|\le 2T$ and 
\begin{equation*} 
K(s)=\sum_nk\bigl(\frac{\log n}{\log P}\bigr)\frac{\Lambda(n)}{\log P}n^{1-s} .
\end{equation*} 
We evaluate $K(s)$  using the Prime Number Theorem in the form 
\begin{equation*} 
\psi(x) =\sum_{n\le x}\Lambda (n)  = x + \Delta (x)\quad {\rm with}\,\, \Delta (x)\ll x(\log x)^{-4} .
\end{equation*} 

We write
\begin{equation*} 
K(s)=(\log P)^{-1}\int k\bigl( \frac{\log x}{\log P}\bigr)x^{1-s}d\bigl(x+\Delta(x)\bigr) = K_1(s) + K_0(s) ,
\end{equation*} 
say, where $K_1(s)$ is the contribution of the main term $x$ and $K_0(x)$ is the contribution of the error term $\Delta (x)$. 

We have two expressions for $K_1(s)$: 
\begin{equation*} 
K_1(s)= \int k(u)P^{(2-s)u}du = Z^{-2}\int k''(u)P^{(2-s)u}du ,
\end{equation*} 
where $Z=(2-s)\log P$. By the inequality 
\begin{equation*} 
\min(A,B/C)\le (A+B)(1+C)^{-1}
\end{equation*} 
we obtain
\begin{equation*} 
|K_1(s)| \le (1+|Z|^2)^{-1}\int\bigl(|k(u)| + |k''(u)|\bigr) P^{(2-\sigma)u}du .
\end{equation*} 
After multiplying this by $P^{{\frac 54}(\sigma -2)}$ as in~\eqref{eq:3.4} we find that $P$ has exponent $(\frac 54 -u)(\sigma -2)$ so we can replace $\sigma = \beta + \beta'$ by its upper bound $1+\min(\beta', \beta) = 1+\beta^{\flat}$, say. We get 
\begin{equation*} 
|K_1(s)| P^{\tfrac 54 (\sigma -2)}\le (1+|Z|^2)^{-1}\int_0^1\bigl(k(\tfrac 54- u) + |k''(\tfrac 54-u)|\bigr) P^{-(1-\beta^{\flat})u}du .
\end{equation*} 

Now we choose the particular crop function
\begin{equation*} 
k(\tfrac 54 - u)  =41(\sin \pi u)^2 = \tfrac{41}{2} (1-\cos 2\pi u)
\end{equation*} 
getting
\begin{equation*} 
\int_0^1 = 41 \int_0^1 \bigl((\sin \pi u)^2 + 2\pi^2|\cos 2\pi u|\bigr)P^{-(1-\beta^{\flat})u} du .
\end{equation*} 
For the last integral we have two easy estimations:
\begin{equation*} 
(1+2\pi^2)/(1-\beta^{\flat})\log P
\end{equation*} 
and
\begin{equation*} 
\int_0^1 \bigl((\sin \pi u)^2 + 2\pi^2|\cos 2\pi u|\bigr)du = \tfrac 12 +4\pi .
\end{equation*} 
The minimum of these two bounds is less than $1+2\pi^2 +\tfrac 12 +4\pi =2(\pi +\tfrac 12)(\pi +\tfrac 32)$ divided by $1+(1-\beta^{\flat})\log P$. Hence we see that the main term contributes to~\eqref{eq:3.4} at most
\begin{equation*} 
1387\sum_{\rho}\sum_{\rho'}|z_{\rho}z_{\rho'}|\bigl(1+(1-\beta^{\flat})\log P\bigr)^{-1}\bigl(1+ (\gamma -\gamma')^2(\log P)^2\bigr)^{-1}.
\end{equation*} 
Finally, applying the inequality $2|z_{\rho}z_{\rho'}| \le  |z_{\rho}|^2 + |z_{\rho'}|^2 $ we get the bound that agrees with the leading term in~\eqref{eq:3.3}. 

To estimate $K_0(s)$ we integrate by parts and proceed as follows:
\begin{equation*} 
\begin{aligned}
K_0(s) & = - (\log P)^{-1}\int\Delta (P^u)d k(u)P^{(1-s)u} \\
& = -\int\bigl( k(u)(1-s) +k'(u)/\log P\bigr)\Delta (P^u)P^{(1-s)u} du \\
& \ll T(\log P)^{-4}\int_{1/4}^{5/4}P^{(2-\sigma)u}du \ll  T(\log P)^{-4}P^{\frac 54 (2-\sigma)} .
\end{aligned}
\end{equation*} 
Hence the contribution to~\eqref{eq:3.4} of the error term $\Delta(x)$ is of order bounded by $T(\log P)^{-4}$ times the quantity 
\begin{equation*} 
\bigl(\sum_{\rho}|z_{\rho}|\bigr)^2\le N\sum_{\rho}|z_{\rho}|^2 .
\end{equation*} 
The result obtained agrees with the final term of~\eqref{eq:3.3}. Our proof of of~\eqref{eq:3.3} is thus complete and hence so, by duality, is that of~\eqref{eq:3.2}. 

\smallskip 

{\bf Corollary 3.2:} {\sl Let $\rho$ run over the zeros of $L(s,\chi)$ with $|\gamma| \le \log q$ and $\beta^*$ denote the maximum of all the $\beta $'s. Let 
$P\ge 
q^6$, $ X\ge P^{5/2}$. 
Then, for any complex numbers $a_p$ supported on $P\le p \le P^{6/5}$
 we have }

\begin{equation}\label{eq:3.5} 
\sum_{\rho}X^{\beta -1}\big|\sum_p a_p p^{-\rho}\big|^2 
\le 2082\bigl(XP^{-5/2}\bigr)^{(\beta^*-1)}\sum_p |a_p|^2p^{-1} . 
\end{equation} 

\smallskip 

{\bf Proof:} Apply the inequality 
$X^{\beta - 1}\le P^{{\frac 52}( \beta -1)}\bigl(XP^{-5/2}\bigr)^{(\beta^*-1)}$. 
Then use~\eqref{eq:3.2} and Corollary A3. Moreover, $NT\le T^2\log qT \ll (\log q)^3$ so the error term in~\eqref{eq:3.2} is absorbed by rounding up 
$1387\cdot 3001/2000 <2082$ .

\section{\bf Prime trios with a character}

Throughout $(a_p)$, $(b_p)$ are sequences of numbers $0\le a_p $, $b_p \le 1$ supported on $P\le p \le P^{6/5}$. Our goal is to estimate the character sum 
\begin{equation}\label{eq:4.1} 
T(X,\chi) =\sum_p\sum_{p_1}\sum_{p_2}\chi(pp_1p_2) a_{p_1}b_{p_2}f(pp_1p_2/X)\log p .
\end{equation} 
Here, $X$ is at our disposal; it is not related to that in~\eqref{eq:2.5}. We assume
\begin{equation}\label{eq:4.2} 
P^{5/2} \le X \le P^4 . 
\end{equation} 
 The crop function $f(y)\ge 0$ is supported on $1\le y \le 2$  with $|f^{(j)}(y)|\le 1$ for $0\le j\le 4$. These conditions imply that the Mellin transform satisfies the bound 
\begin{equation*} 
\tilde f(s)=\int f(y)y^{s-1}dy \ll(1+|s|)^{-4} 
\end{equation*} 
and that 
\begin{equation*} 
|\tilde f(s)|\le |\tilde f(1)|= \int f(y)dy =\hat f(0)
\end{equation*} 
when $0\le {\rm Re}\, s \le 1 $.

For the principal character $\chi = \chi_0$ we have 
\begin{equation}\label{eq:4.3} 
T(X;1) =\bigl(\hat f (0) +\tau\bigr)X\ab \bb
\end{equation} 
where 
\begin{equation}\label{eq:4.4} 
\ab = \sum_pa_p p^{-1} , \quad \bb = \sum_pb_p p^{-1} .
\end{equation} 
Note that $ \ab$, $\bb \le \log {\frac65} +\tau$ and that equality holds if 
$a_p = b_p =1$. In several places we shall use the estimates 
$\frac{1}{31} <  (\log {\frac65} )^2 <\frac{1}{30}$, 
so that we can drop the term with $\tau =O(1/\log P)$. 

\smallskip 

{\bf Lemma 4.1:} {\sl For $\chi \neq \chi_0$ we have 
\begin{equation}\label{eq:4.5} 
|T(X;\chi)| \le 380 \hat f (0)X^{\beta^*} P^{\frac 52 (1-\beta^*)}
\end{equation} 
where }
\begin{equation}\label{eq:4.6} 
\beta^* =\max\{\beta; \, \rho = \beta +i\gamma,  \, L(\rho, \chi)=0,  \, |\gamma| \le \log q\} .
\end{equation} 

\smallskip 

{\bf Proof:} We use the approximate explicit formula 
\begin{equation*} 
\sum_p\chi(p)f(p/Y)\log p = -\sum_{|\gamma|\le T}\tilde f (\rho)Y^{\rho} + O\bigl(Y^{\beta^*}(\log q)^{-2}\bigr)
\end{equation*} 
with $T=\log q$ and $Y=X/p_1p_2$, getting 
\begin{equation*} 
\begin{aligned}
T(X;\chi) =  -\sum_{|\gamma|\le T}\tilde f (\rho)X^{\rho} & 
\bigl(\sum_p\chi(p)a_pp^{-\rho}\bigr)\bigl(\sum_p\chi(p)b_pp^{-\rho}\bigr) \\
& + O\bigl(X^{\beta^*}P^{\frac{12}{5}(1-\beta^*)}(\log q)^{-2}\bigr) .
\end{aligned}
\end{equation*} 
Hence~\eqref{eq:4.5} follows from~\eqref{eq:3.5} by Cauchy's inequality and rounding up $2082\cdot \log \tfrac 65 < 380$.

\smallskip 

Next, for a real character $\chi\neq \chi_0$ we are going to approximate $T(X;\chi)$ by $-T(X;\chi_0)$. Put 
\begin{equation}\label{eq:4.7} 
\lambda(p) = (1 * \chi)(p). 
\end{equation}
This is unrelated to our constant $\lambda >1$  introduced earlier. Note that 
$1+\chi(pp_1p_2) \le \lambda(p) +\lambda(p_1)+\lambda(p_2)$. Hence, 
$T(X;\chi) +T(X;\chi_0)$ is bounded by
\begin{equation}\label{eq:4.8} 
\sum_p\sum_{p_1}\sum_{p_2}\bigl(\lambda(p) +\lambda(p_1)+\lambda(p_2)\bigr) a_{p_1}b_{p_2}f(pp_1p_2/X)\log p .
\end{equation} 
The terms with $\lambda(p_1)$  contribute to~\eqref{eq:4.8} exactly~\eqref{eq:4.3}, but with the coefficient $a_{p_1}$ replaced by  $\lambda(p_1)a_{p_1}$, similarly with $\lambda(p_2)$. 

To estimate the contribution of the terms with $\lambda(p)$, we first execute the summation over $p_2$ in~\eqref{eq:4.8}. Note that $\log p < (2+\tau)\log p_2$ because $p_2\ge P$ and $p\le 2XP^{-2} \le 2P^2$. Moreover, if $P>q^{20}$ we have 
$p\ge XP^{-12/5}\ge P^{1/10}> q^2$. Hence, the contribution of $\lambda (p)$ in~\eqref{eq:4.8} is bounded by
\begin{equation*} 
(2+\tau)\sum_{q^2<p \le2P^2}\lambda(p)\sum_{p_1}\sum_{p_2}a_{p_1}f(pp_1p_2/X)\log p_2\le (2+\tau)\hat f (0)X\ab\delta(2P^2)
\end{equation*} 
where
\begin{equation}\label{eq:4.9} 
\delta(z)=\sum_{q^2<p\le z}\lambda (p) p^{-1} .
\end{equation} 
Interchanging $a_p$ with $b_p$ we can replace $\ab$ above by $(\ab \bb)^{\frac12}$. Adding up the three contributions we see that $T(X;1)+T(X;\chi)$ is bounded by
\begin{equation*} 
(\hat f (0) +\tau) X\bigl\{ \ab_{\lambda} \bb + \ab \bb_{\lambda} + 2(\ab \bb)^{\frac12} \delta(2P^2)\bigr\} 
\end{equation*} 
where 
\begin{equation*} 
\ab_{\lambda}=\sum_p\lambda(p)a_pp^{-1}\le \delta(P^{6/5})\le \tfrac{12}{5}(1-\beta)\log P
\end{equation*} 
and similarly $\bb_{\lambda}\le \frac{12}{5}(1-\beta)\log P$, see (A5). Hence, we conclude: 

\smallskip 

{\bf Lemma 4.2:} {\sl  For $\chi$ real, $\chi \neq \chi_0$, $P\ge q^{20}$ and 
$P^{5/2} \le X\le P^4$ we have
\begin{equation}\label{eq:4.10} 
0\le T(X;1) + T(X; \chi) \le \tfrac{12}{5} \hat f (0)X(1-\beta)\log P
\end{equation} 
where $\beta$ is any real zero of $L(s, \chi)$. }
 
\smallskip

{\bf Remarks:} For example, $\beta <0$ could be a trivial zero, but then~\eqref{eq:4.10} would not be interesting. The constant $12/5$ in~\eqref{eq:4.10} comes by a rounding up of 
$2\cdot 2 \cdot (\tfrac 65 +2) \log \tfrac 65 = 2.3337\ldots$. Combining~\eqref{eq:4.10} with~\eqref{eq:4.3}, one obtains another bound for $T(X; \chi) $ for $\chi$ real, in addition to~\eqref{eq:4.5}. This bound is better than that given in~\eqref{eq:4.5} if $L(s,\chi)$ has a real zero $\beta$ very close to $1$.

\section{\bf Prime quintets in arithmetic progressions}

Throughout, we shall assume the conditions of Section 4. In addition to $(a_p)$, $(b_p)$ supported on $P\le p \le P^{6/5}$ we consider two sequences $(c_p)$, $(d_p)$ of numbers $0\le c_p, d_p \le 1$ which are supported on short segments $C<p\le \lambda C$, $D< p\le \lambda D$ respectively; see~\eqref{eq:2.32}. Note that the choice 
\begin{equation}\label{eq:5.1} 
X=x/CD
\end{equation} 
satisfies $X\le xP^{-2}=P^4=x^{2/3}$ and $X>xP^{-12/5} =P^{18/5}=x^{3/5}$. 

Our goal is to estimate the sums
\begin{equation}\label{eq:5.2} 
Q(x;q,a)= \sum_p\ssssum_{pp_1p_2p_3p_4\equiv a (\mod q)}a_{p_1}b_{p_2}c_{p_3}d_{p_4}f(pp_1p_2/X) \log p
\end{equation} 
with the same crop function $f(y)$ as in Section 4. We have 
\begin{equation}\label{eq:5.3} 
Q(x;q,a)= \frac{1}{\vp(q)}\sum_{\chi(\mod q)}\bar\chi (a) T(X;\chi)\cC(\chi)\cD(\chi)
\end{equation} 
where $T(X;\chi)$ was defined in~\eqref{eq:4.1} and 
\begin{equation}\label{eq:5.4} 
\cC(\chi)=\sum_p\chi(p)c_p , \quad \quad \cD(\chi)=\sum_p\chi(p)d_p .
\end{equation} 

Similarly, we define $\cC(1), \cC(\lambda), \cD(1), \cD(\lambda)$ which will appear later. We also introduce 
\begin{equation}\label{eq:5.5} 
\omega(\cC)=\sum_{C<p\le\lambda C}p^{-1} , \quad \quad \frak S (C)=\sum_{C<p\le\lambda C}\lambda(p)p^{-1} . 
\end{equation} 
Similarly, we define $\omega(\cD), \frak S (D)$. 

\smallskip

Note that $\cC(1)\le\lambda C\omega(C)$,  $\cD(1)\le\lambda D\omega(D)$. Moreover, if $c_p=d_p=1$, then we have lower bounds $\cC(1)\ge C\omega(C)$, $\cD(1)\ge D\omega(D)$.  

\smallskip

We have treated $T(X;\chi)$ in Section 4; see~\eqref{eq:4.3},~\eqref{eq:4.5} and~\eqref{eq:4.10}. The last estimate~\eqref{eq:4.10} applies to the real character $\chi\neq \chi_0$.  We compared $T(X;\chi)$ to $-T(X;1)$ apart from a small quantity of order $X(1-\beta) \log P$.  Here, we make a similar comparison of $\cC(\chi)\cD(\chi)$ to $\cC(1)\cD(1)$. For every $p, p'$ we have 
$0\le 1- \chi(pp') \le \lambda(p) + \lambda(p')$, so
\begin{equation}\label{eq:5.6} 
\begin{aligned}
0\le \cC(1)\cD(1)-\cC(\chi)\cD(\chi) & \le \cC(\lambda)\cD(1)+\cC(1)\cD(\lambda)\\ & \le \lambda^2CD\bigl(\frak S (C)\omega(D)+ \omega (C)\frak S (D)\bigr) .
\end{aligned}
\end{equation} 

Now, we return to the formula~\eqref{eq:5.3}. The principal character $\chi_0 (\mod  q)$ yields the main contribution (see~\eqref{eq:4.3})
\begin{equation}\label{eq:5.7} 
Q_0(X)= \frac{1}{\vp(q)}T(X;1)\cC(1)\cD(1) = \bigl({\hat f} (0) +\tau\bigr)\frac{X}{\vp(q)}\ab\bb\cC(1)\cD(1) .
\end{equation} 
There may be a real character $\chi_1\neq \chi_0$ for which its contribution 
\begin{equation}\label{eq:5.8} 
Q_1(X)= \frac{\chi_1(a)}{\vp (q)}T(X;\chi_1)\cC(\chi_1)\cD(\chi_1) 
\end{equation} 
requires special attention. The other characters contribute to~\eqref{eq:5.3} in total at most
\begin{equation}\label{eq:5.9} 
Q^*(X)= \frac{1}{\vp (q)}T_{\max}(X)
\bigl(\sum_{\chi}|\sum_p\chi(p)c_p|^2\bigr)^{\frac 12}
\bigl(\sum_{\chi}|\sum_p\chi(p)d_p|^2\bigr)^{\frac 12}
\end{equation} 
where $T_{\max}(X)$ denotes the maximum of $|T(X,\chi)|$ over all $\chi$
other than $\chi_0$ or $\chi_1$ (if the latter exists). For these non-exceptional characters we have 
\begin{equation}\label{eq:5.10}
\beta^*\le 1-\eta/\log q
\end{equation} 
by Lemma A1, where $\eta$ is an absolute constant to be specified later. Hence, Lemma 4.1, together with the bounds $P^{5/2}X^{-1}  < x^{5/12 - 3/5} < x^{-1/6}$ yield
\begin{equation}\label{eq:5.11}
|T_{\max}(X)|\le 380{\hat f}(0) X x^{-\eta/ 6\log q} .
\end{equation} 
Next, we apply Lemma A5 with $C\ge P$, so $C/q\ge P/q\ge P^{19/20}$ and~\eqref{eq:5.9} yields
\begin{equation}\label{eq:5.12}
\vp (q)Q^*(X)\le 801 {\hat f}(0)\bigl(\frac{\log\lambda}{\log P}\bigr)^2 x^{1-\eta/ 6\log q} .
\end{equation} 
Here, 801 is obtained after increasing the number $380\cdot 2\cdot 20/19 =800$, which allows us to replace $\lambda -1$ by $\log \lambda$.
If the exceptional character $\chi_1$ exists, which means that $L(s, \chi_1)$ has a real zero
\begin{equation}\label{eq:5.13}
  \beta_1 > 1- \eta /\log q ,
\end{equation} 
we evaluate its contribution by comparing $\vp (q)Q_1(X)$ to
\begin{equation}\label{eq:5.14}
 \vp (q)Q_{11}(X)= -\chi_1(a)T(X;1)\cC(1)\cD(1) .
\end{equation} 
We get 
\begin{equation}\label{eq:5.15}
  \begin{aligned}
    \vp (q)|Q_{1}(X)-Q_{11}(X)| & \le \bigl(T(X;1)+T(X;\chi_1)\bigr)\cC(1)\cD(1) \\ & \quad \quad + T(X;1)\bigl(\cC(1)\cD(1)-\cC(\chi_1)\cD(\chi_1) \bigr)\\
    & \le\tfrac 25 {\hat f} (0)x \omega (C)\omega (D)(1-\beta_1) \log x \\ &
 \quad \quad   +\tfrac{1}{30} {\hat f} (0)x \bigl(\omega (C) {\frak S}(D)+ \omega (D){\frak S}(C)\bigr) \\
    & \le \tfrac{1}{30} {\hat f} (0)x\Omega(C,D)
    \end{aligned}
\end{equation}
by~\eqref{eq:4.10},~\eqref{eq:5.6} and~\eqref{eq:4.3}, in which $\ab \bb <\frac{1}{30}$. Here
\begin{equation}\label{eq:5.16}
  \Omega(C,D) = 12  \omega(C)\omega(D)(1-\beta_1) \log x + \omega (C) {\frak S}(D)+ \omega (D){\frak S}(C) .
  \end{equation} 

At some point (see~\eqref{eq:6.6}) we shall have to sum over $C$, $D$ in the segment~\eqref{eq:2.33}. Note that
\begin{equation*}
\sum_C\omega(C)=\sum_{P<p\le P^{6/5}}p^{-1} =\log \frac 65 +\tau ,
  \end{equation*} 

\begin{equation*}
\sum_C{\frak S}(C)=\sum_{P<p\le P^{6/5}}\lambda (p)p^{-1} \le \delta(x^{1/5}) 
  \end{equation*} 
with $ \delta(x^{1/5}) \le \tfrac 25 (1-\beta_1) \log x$ by Lemma A4. The same estimates hold for the sums over $D$. Hence,
 \begin{equation}\label{eq:5.17}
  \sum_C\sum_D\Omega(C,D) \le \frac{11}{20}(1-\beta_1) \log x ,
  \end{equation} 
 where $11/20$ comes from rounding up $12(\log \tfrac65)^2 + \tfrac 45 \log \tfrac65$ .

 If $\chi_1$ is exceptional we can still use Lemma 4.1 which gives
 \begin{equation*}
   |T(X;\chi_1)| \le 380 {\hat f}(0)Xx^{-(1-\beta_1)/6}
 \end{equation*}
 in place of~\eqref{eq:5.11} . Hence~\eqref{eq:5.8} gives the bound 
\begin{equation}\label{eq:5.18}
  \varphi(q)Q_1(X) \le 380 {\hat f}(0)\bigl(\frac{\log \lambda}{\log P}\bigr)^2
  x^{1-(1-\beta_1)/6} 
   \end{equation} 
in place of~\eqref{eq:5.12}, which is still useful if the exceptional zero $\beta_1$ is not very close to 1.

\section{\bf Three estimations of $Q(\cA)$}

Recall that $Q(\cA)$ is the sum over prime quintets $n=pp_1p_2p_3p_4$ in the aritmetic progression $n\equiv a (\mod q)$ counted with the smooth weight $f(n/x)$; see~\eqref{eq:2.31},~\eqref{eq:2.37},~\eqref{eq:2.38}. Now we use the results of the previous two sections to estimate $Q(\cA)$ in terms of the zeros of $L(s,\chi)$. We start from~\eqref{eq:5.2} and~\eqref{eq:5.3} which we have estimated by:
\begin{equation}\label{eq:6.1}
Q(x;q,a)\ge Q_0(X) - Q^*(X)- Q_1(X);
   \end{equation} 
see~\eqref{eq:5.7},~\eqref{eq:5.9},~\eqref{eq:5.8}, with $X=x/CD$. 
The first part $Q_0(X)$ comes from the principal character. The second part, which comes from those non-principal characters other than the exceptional one, satisfies the bound~\eqref{eq:5.12} which is sufficiently strong because $\eta$ is an absolute constant, which can be taken to be not very small. The last part $Q_1(X)$ comes from the exceptional character $\chi_1$ and it satisfies the bound~\eqref{eq:5.18}, which is good enough only if the exceptional zero $\beta_1$ is extremely close to 1. In case it is extremely close, then $Q_1(X)$ is comparable to $Q_{11}(X)$; see~\eqref{eq:5.15}. Therefore, we need a second option, in addition to the direct one~\eqref{eq:5.18} for the estimatiom of $\vp(q)Q_1(X)$. This one is provided by
\begin{equation}\label{eq:6.2}
\vp(q)Q_1(X) \le -\chi_1(a)T(X;1)\cC(1)\cD(1) + \tfrac{1}{30}\hat f (0)x\Omega(C,D); 
   \end{equation} 
see~\eqref{eq:5.14},~\eqref{eq:5.15}. Hence, if $\chi_1$ exists then, using~\eqref{eq:5.7} for $Q_0(X)$,
\begin{equation}\label{eq:6.3}
\vp(q)Q(x;q,a) ) \ge \hat f (0)x\cW(C,D) 
   \end{equation} 
where we have two choices:
\begin{equation}\label{eq:6.4}
\cW(C,D) =\frac{1}{31} 
\omega(C)\omega(D)-801\bigl(\frac{\log \lambda}{\log P}\bigr)^2x^{-\eta/6\log q}
-380\bigl(\frac{\log \lambda}{\log P}\bigr)^2x^{(1-\beta_1)/6}
   \end{equation} 
and
\begin{equation}\label{eq:6.5}
\cW(C,D) =\frac{1}{31} (1-\chi_1(a))\omega(C)\omega(D)-801\bigl(\frac{\log \lambda}{\log P}\bigr)^2x^{-\eta/6\log q}
-\frac{1}{30}\Omega(C,D) .
   \end{equation} 

In these, we obtained the fraction $1/31$ on rounding down $(\log 6/5)^2>1/31$. If $\chi_1$ does not exist then $\cW(C,D)$ is given by~\eqref{eq:6.4} without the last term. Recall that the number of pairs $C$, $D$ is $(\log P)^2(5\log\lambda)^{-2}$. 

\smallskip

Introducing the factor $\log p\le\log(2x^{1/3})$ into~\eqref{eq:2.36} we get
\begin{equation}\label{eq:6.6}
\vp(q)Q_f(C,D)\log(2x^{1/3})\ge \vp(q)Q(x;q,a)
\ge \hat f (0)x\cW(C,D) .
   \end{equation} 
Note we have omitted the negative contribution of $Q_h(C,D)$ in~\eqref{eq:2.37} because it is very small for $\lambda$ close to 1 and so is absorbed by the margin we obtained in rounding up the constants in~\eqref{eq:6.4},~\eqref{eq:6.5}. Next, summing over $C$ and $D$ in~\eqref{eq:6.4},~\eqref{eq:6.5} and using ~\eqref{eq:5.17}, we obtain respectively two lower bounds
\begin{equation*}
\frac{1}{961}-\frac{801}{25}x^{-\eta/6\log q}-\frac{380}{25}x^{-(1-\beta_1)/6}
   \end{equation*} 
and
\begin{equation*}
\frac{1}{961} (1-\chi_1(a))-\frac{801}{25}x^{-\eta/\log q}-\frac{11}{600}
(1-\beta_1)\log x .
   \end{equation*} 
We simplify these bounds by assuming
\begin{equation}\label{eq:6.7}
x\ge q^{80/\eta} 
   \end{equation} 
which makes the middle term smaller than one twentieth of the first term $1/961$ (check that $801/25e^{80/6}< 1/19270<1/961\cdot20$). Hence, we are left with two lower bounds
\begin{equation}\label{eq:6.8}
321\vp(q)Q(\cA)\log x \ge\hat f(0)x\{\tfrac{19}{20}-14608x^{(1-\beta_1)/6}\} 
   \end{equation} 
and
\begin{equation}\label{eq:6.9}
321\vp(q)Q(\cA)\log x \ge\hat f(0)x\{\tfrac{19}{20}-\chi_1(a)- 18(1-\beta_1)\log x\} .
   \end{equation} 
If the exceptional character  $\chi_1$ does not exist, then~\eqref{eq:6.8} holds without the negative exceptional term. We have proved the following result.

\smallskip

{\bf Lemma 6.1:} {\sl Suppose that all the zeros $\rho=\beta +i\gamma$ of every $L(s,\chi)$ satisfy 
\begin{equation}\label{eq:6.10}
\beta\le 1-\eta/\log q(|\gamma|+1) ,
   \end{equation} 
except possibly for one simple real zero $\beta_1$ of $L(s,\chi_1)$ with a real character $\chi_1$. Let $x\ge q^{80/\eta}$. Then we have
\begin{equation}\label{eq:6.11}
\vp(q)Q(\cA)\ge \hat f(0)x/350\log x ,
   \end{equation} 
subject to any one of the following conditions:
\begin{equation}\label{eq:6.12}
\chi_1 \,\, {\sl does \,\, not \,\, exist} ,
   \end{equation} 
\begin{equation}\label{eq:6.13}
\chi_1 \,\, {\sl exists\,\, and}\,\, x\ge e^{80/(1-\beta_1)}
   \end{equation} 
\begin{equation}\label{eq:6.14}
\chi_1 \,\, {\sl exists}, \,\, \chi_1(a)=-1\quad {\sl and}\quad x\le e^{1/18(1-\beta_1)} .
   \end{equation} }

\medskip

{\bf Proof:} The exceptional term in~\eqref{eq:6.8}, subject to~\eqref{eq:6.13}, is smaller than $14608e^{-80/6} < 1/42$. Check that $19/20 - 1/42 > 321/350$. This shows that~\eqref{eq:6.11} holds subject to~\eqref{eq:6.13}. The other two cases are clear. 

\medskip 

{\bf Remark.} The true expected value of $\vp(q)Q(\cA)$ was computed just after~\eqref{eq:2.31} so our lower bound~\eqref{eq:6.11} is about one third that size. 

\section{\bf Primes in an arithmetic progression}

Finally, we proceed to the estimation of the sum over primes~\eqref{eq:2.28}. To this end, we appeal to the sieve inequality~\eqref{eq:2.29} which transfers the task to that for prime quintets apart from a small piece
\begin{equation}\label{eq:7.1}
\hat f (0)\bigl(\log \frac{1}{2\theta -1}\bigr)\frac{2x}{\theta \vp(q)\log x} 
   \end{equation} 
accounting for prime duos, which needs to be subtracted. Recall that $y=x/q(\log x)^3= x^{\theta}$. We make the contribution in~\eqref{eq:7.1} smaller than $\tfrac{1}{25}Q(\cA)$ by choosing $x\ge q^M$ with a sufficiently large number $M$. In view of the lower bound~\eqref{eq:6.11}, any $M$ satisfying 
\begin{equation*}
\frac{2M}{M-1}\log\frac{M}{M-3} <\frac{1}{25\cdot 350}
   \end{equation*} 
will suffice. We check that $M= 52600$ is good enough. Applying~\eqref{eq:2.29} and~\eqref{eq:6.11} and checking that 1/24 - 1/25 = 350/210000, we find:

\medskip

{\bf Lemma 7.1:} {\sl If $x\ge q^{52600}$ and one of the conditions~\eqref{eq:6.12},~\eqref{eq:6.13},~\eqref{eq:6.14} holds, then
\begin{equation}\label{eq:7.2}
\vp(q)\pi_f(x;q,a) \ge \hat f(0)x/210000\log x .
   \end{equation} }

\smallskip 

In addition to~\eqref{eq:7.2}, we have the lower bound
\begin{equation}\label{eq:7.3}
\vp(q)\pi_f(x;q,a) \ge L(1,\chi_1)V(\chi_1)x/168 
   \end{equation} 
if $\chi_1(a) =1$ and $q^{43}\le x \le e^{1/4(1-\beta_1)}$.
This bound is derived in [FI3] by quite different arguments using Selberg's lower bound sieve. Combining~\eqref{eq:7.2} and~\eqref{eq:7.3} we establish: 

\smallskip

{\bf Theorem 7.2:} 
{\sl Let $q$ be sufficiently large and $(a,q)=1$. We have
\begin{equation}\label{eq:7.4}
p_{\min}(q,a) \le q^L \,\, {\rm with} \,\, L = 75744000. 
   \end{equation} }

\smallskip

{\bf Proof:} If the exceptional zero does not exist then~\eqref{eq:7.2} produces~\eqref{eq:7.4} with $L=M=52600$. If the exceptional zero $\beta_1$ exists and satisfies $(1-\beta_1)\log q\ge 1/18\cdot 52600$, then~\eqref{eq:7.2} subject to~\eqref{eq:6.13} together with~\eqref{eq:7.3} subject to $x\le e^{1/4(1-\beta_1)}$ yield~\eqref{eq:7.4} with $L=M=52600$. In the opposite case $(1-\beta_1)\log q < 1/18\cdot 52600$, we use~\eqref{eq:7.2} subject to~\eqref{eq:6.13}, yielding~\eqref{eq:7.4} with $L=80\cdot 18\cdot 52600 = 75744000$.

\medskip

{\bf Remarks:} We assumed that $q$ is sufficiently large and $q^{80/\eta}\le q^{52600}$,that is the non-exceptional zeros satisfy~\eqref{eq:6.10} with $\eta \ge 1/675.5$ On the other hand, the old result of R.J. Miech [Mi] provides the much larger constant $\eta = 1/20$. The current record, given specifically for $|t| \le 1$, is due to Xylouris [Xy] and stands at $\eta = .440$.

\medskip 

\quad\quad\quad\quad\quad\quad\quad\quad\quad\quad\quad\quad {\bf Appendix}

\medskip 

In this section we select some classical results about the zeros of Dirichlet $L$-functions and character sums over prime numbers.

\medskip

{\bf Lemma A1:} For every $\chi (\mod q)$, all zeros $\rho =\beta +i\gamma$ of $L(s,\chi)$ satisfy 
\begin{equation*}%\label{eq:8.1}
\beta\le 1-\eta/\log q(1+|\gamma|) ,
   \end{equation*} 
except possibly for one real simple zero $\beta_1$ of $L(s,\chi_1)$ with one real character $\chi_1(\mod q)$. Here, $\eta$ is an absolute positive constant. 

\medskip

We need the above inequality %~\eqref{eq:8.1} 
with $\eta=1/657$. 
We do not need the log-free zero density estimation nor the Deuring-Heilbronn repulsion property. In place of these, we quickly prove the following inequality.

\medskip
 
{\bf Lemma A2:} For $s=\sigma + it$ with $\sigma >1$ we let
\begin{equation*}%\label{eq:8.2}
V(s,\chi) =\sum_{\rho}\bigl(1+\frac{1-\beta}{\sigma - 1}\bigr)^{-1}
\bigl(1+\bigl(\frac{\gamma -t}{\sigma - 1}\bigr)^2\bigr)^{-1} ,
\end{equation*} 
where $\rho =\beta +i\gamma$ runs over the zeros of $L(s,\chi)$ with $\beta>0$. If $\chi \neq \chi_0$ we have
\begin{equation*}%\label{eq:8.3}
V(s,\chi) \le  1 + \frac{\sigma - 1}{2} \log cq|s| 
\end{equation*}
where $c>1$ is an absolute constant.

 \medskip 

{\bf Proof:} We can assume $\chi$ is primitive. Then (see e.g. (5.24) of [IK])
\begin{equation*}
\frac{L'}{L}(s,\chi) = -\frac 12 \log \frac{q}{\pi} -\frac 12 \frac{\Gamma'}{\Gamma}(\frac{s_{\chi}}{2}) +B(\chi)+\sum_{\rho}\bigl(\frac{1}{s-\rho}+ \frac{1}{\rho}\bigr)
\end{equation*}
where $s_{\chi}= s +(1+\chi(-1))/2$, $\Gamma'(s)/\Gamma(s)=\log s +O(1)$ and 
\begin{equation*}
{\rm Re}\, B(\chi) = -\sum_{\rho}{\rm Re}\,\frac{1}{\rho} .
\end{equation*}
Hence,
\begin{equation*}
\sum_{\rho}{\rm Re}\,\frac{1}{s-\rho} = \frac12 \log q|s| +{\rm Re}\,\frac{L'}{L}(s,\chi) +O(1) ,
\end{equation*}
\begin{equation*}
\bigl|\frac{L'}{L}(s,\chi)\bigr|\le -\frac{\zeta'}{\zeta}(\sigma) =\frac{1}{\sigma -1} +O(1) .  
\end{equation*}
Now Lemma A2  %~\eqref{eq:8.3} 
  follows easily from the inequality
\begin{equation*}
{\rm Re}\,\frac{1}{s-\rho}=\frac{\sigma -\beta}{(\sigma -\beta)^2 +(\gamma -t)^2}\ge 1/(\sigma-1)\bigl(1+\frac{1-\beta}{\sigma -1}\bigr)\bigl(1+\bigl(\frac{\gamma -t}{\sigma -1}\bigr)^2\bigr) .
\end{equation*}

\medskip

{\bf Corollary A3:} For $q$ sufficiently large and $|t| \le \log q$ we have
\begin{equation*}%\label{eq:8.4}
\sum_{\rho}\bigl(1+(1-\beta)\log q\bigr)^{-1}\bigl(1+(\gamma - t)^2(\log q)^2\bigr)^{-1} \le \frac 32 +\frac{1}{2000} .
\end{equation*}

\medskip 

{\bf Lemma A4:} Let $\chi (\mod q)$ be a real character and $\beta$ a real zero of $L(s, \chi)$. Then, for $x\ge q^2$, $q$ sufficiently large, we have
\begin{equation*}%\label{eq:8.5}
\sum_{q^2\le p \le x}\bigl(1+\chi(p)\bigr)p^{-1} \le 2(1-\beta) \log x .
\end{equation*}

\smallskip

{\bf Proof:} See the text between Proposition 24.1 and Corollary 24.2 of [FI2].

\medskip 

{\bf Lemma A5:} If $C\ge q^2$ and $|c_p|\le 1$ for $C< p \le \lambda C$, then 
\begin{equation*}%\label{eq:8.6}
\sum_{\chi (\mod q)}\bigl|\sum_{C < p \le \lambda C}\chi (p) c_p\bigr|^2 \le (2+\tau) (\lambda -1)^2C^2 (\log C)^{-1}(\log C/q)^{-1}
\end{equation*}
where $\lambda >1$ is a fixed number and $\tau \ll (\log C)^{-1}$ . 

\smallskip 

{\bf Proof:} The left side is estimated by the Brun-Titchmarsh theorem as follows:
\begin{equation*}
\begin{aligned}
\vp (q)\ssum_{p' \equiv p (\mod q)} c_{p'}{\bar c_p} & \le \vp(q)\sum_{C < p \le \lambda C}
\bigl(\pi (\lambda C; q, a) - \pi (C; q, a)\bigr) \\
& \le (2+\tau) (\lambda -1 )C (\log C/q)^{-1} \bigl(\pi (\lambda C) - \pi (C)\bigr)  \\
& \le (2+\tau) (\lambda -1 )^2C^2(\log C/q)^{-1}(\log C)^{-1} .
\end{aligned}
\end{equation*}

\medskip

{\bf Remark:} Our treatment of character sums over prime quintets resembles a multiplicative version of the circle method for the ternary Goldbach problem. The prime trios alone would perform this part, but the extra two prime factors are needed to avoid the log-free zero density bound~\eqref{eq:1.3} and the repulsion property of the possible exceptional zero. There is a possibility to arrange a combinatorial identity of sieve type by means of which we could work with two primes and one almost prime rather than with prime quintets.

\medskip 
Department of Mathematics, University of Toronto

Toronto, Ontario M5S 2E4, Canada  \quad (frdlndr@math.toronto.edu)

\medskip

Department of Mathematics, Rutgers University

Piscataway, NJ 08903, USA  \quad  (iwaniec@comcast.net)

\end{document}